# Solution Representations for a Wave Equation with Weak Dissipation

Jens Wirth[*]

October 27, 2018

## Abstract

We consider the Cauchy problem for the weakly dissipative wave equation

$$\Box v + \frac{\mu}{1+t} v_t = 0, \qquad x \in \mathbb{R}^n, \quad t \geq 0,$$

parameterized by $\mu > 0$, and prove a representation theorem for its solutions using the theory of special functions.

This representation is used to obtain $L_p$–$L_q$ estimates for the solution and for the energy operator corresponding to this Cauchy problem.

Especially for the $L_2$ energy estimate we determine the part of the phase space which is responsible for the decay rate. It will be shown that the situation depends strongly on the value of $\mu$ and that $\mu = 2$ is critical.

## 1 Introduction

We are interested in properties of solutions to the Cauchy problem

$$\Box v + \frac{\mu}{1+t} v_t = 0, \quad v(0,\cdot) = v_1, \quad v_t(0,\cdot) = v_2 \qquad (1.1)$$

to data $v_1, v_2 \in \mathcal{S}(\mathbb{R}^n)$.

---

[*]Institute of Applied Mathematics I, TU Bergakademie Freiberg, 09596 Freiberg, Germany, email: wirth@math.tu-freiberg.de

**I.** On the one hand this problem is the model case of weak dissipation. The problem is investigated by this reason, but not well understood. This shall be underlined by one result. If we consider the usual hyperbolic energy of a solution

$$E(v;t) = \frac{1}{2}\int_{\mathbb{R}^n}\left(v_t^2 + (\nabla v)^2\right)\mathrm{d}x, \tag{1.2}$$

it is known from [Mat77], [Ues79] that

$$E(v;t) = O(t^{-\alpha}), \qquad \alpha = \min\{2,\mu\}. \tag{1.3}$$

The method of proof are weighted energy inequalities. They provide a powerful tool of proving such a-priori estimates without knowing details on the solution representation.

So far it is not motivated why this change in the decay order for $\mu = 2$ appears. We will give an alternative proof of such a-priori estimates using an explicit representation of the solution operator by means of special functions. Later on we will see that different techniques are necessary to understand small and large values of the parameter $\mu$.

**II.** On the other hand this equation is related to wave equations with increasing in time propagation speed. If we are interested in solutions to the Cauchy problem for the operator

$$\partial_\tau^2 - \lambda^2(\tau)\Delta, \tag{1.4}$$

investigated e.g. in [RY98], we can apply a change of coordinates to reduce it to an equation with constant speed of propagation and dissipative term. We introduce the new time variable

$$t := \int_{\tau_0}^\tau \lambda(s)\mathrm{d}s,$$

such that $\partial_\tau = \lambda(\tau)\partial_t$ and $\partial_\tau^2 = \lambda^2(\tau)\partial_t^2 + \lambda'(\tau)\partial_t$. Thus in the variable $t$ the operator rewrites as

$$\lambda^2(\tau)\left[\Box + \frac{\lambda'(\tau)}{\lambda^2(\tau)}\partial_t\right]. \tag{1.5}$$

Especially for the coefficient $\lambda(t) = (1+t)^\ell$, investigated by M. Reissig in [Rei97], we get with $\tau_0 = (1+\ell)^{\frac{1}{1+\ell}} - 1$ the relation $t = 1/(\ell+1)(1+\tau)^{\ell+1} - 1$ and therefore the equivalent problem

$$\Box + \frac{\ell}{\ell+1}\frac{1}{1+t}\partial t,$$

i.e. (1.1) with $\mu = \ell/(\ell+1) \in (0,1)$.

Similarly for $\lambda(t) = e^t$ we get the relation to our equation with parameter $\mu = 1$. This Cauchy problem was recently studied by A. Galstian in [Gal00]. This paper extents the results of [Rei97] and [Gal00] and provides in the case of the latter one even sharper estimates.

In all the papers [Rei97], [Gal00] and [TT80] the Cauchy problem was reduced to the confluent hypergeometric equation by change of coordinates in frequency space. We will follow this idea.

**III.** It should be remarked that the related singular problem in one space dimension
$$v_{tt} + \frac{\mu}{t} v_t - v_{xx} = 0$$
is related to the free wave equation for radial data (understanding $x$ as time variable and $t$ as radius). This follows directly from the representation of the Laplacion $\Delta$ in polar coordinates
$$\Delta = \partial_r^2 + \frac{n-1}{r} \partial_r + \frac{1}{r^2} \Delta_S,$$
where $\Delta_S$ is the Laplace-Beltrami operator on the unit sphere $S^{n-1}$ of $\mathbb{R}^n$.

Thus this special problem generalizes the Euler-Darboux equation to non-integral $\mu$.

In [Wei54] a solution representation for the multi-dimensional analogon
$$v_{tt} + \frac{k}{t} v_t - \Delta v = 0$$
was obtained using recursion formulae and a generalized method of descent. The difference to our representation is that we use Fourier multiplier and work in frequency space while A. Weinstein used integral formulae in physical space.

**IV.** In the scope of this paper are two main points.

At first we give an explicit representation of the solutions to the Cauchy problem (1.1) in terms of Bessel functions. This will be done in Section 2.

Afterwards we use this representation to describe the Fourier multiplier corresponding to
$$(v_1, v_2) \mapsto (\nabla v, v_t)$$
in order to derive energy estimates in the (strong) sense of a norm estimate for this operator. The same is done for the solution operator of the Cauchy problem itself. This is the main contents of Section 3.

In Section 4 we give a further a-priori estimate including a non-trivial source term.

## 2 Multiplier Representation

**Reduction to Bessel's equation.** At first we construct the fundamental solution of the corresponding ordinary differential equation in the Fourier image. Let $\hat{v}(t,\xi) = \mathcal{F}_{x\to\xi}[v]$,

$$\hat{v}(t,\xi) = (2\pi)^{-\frac{n}{2}} \int_{\mathbb{R}^n} e^{-ix\cdot\xi} v(t,x) \mathrm{d}x. \tag{2.1}$$

Then $\hat{v}$ satisfies the ordinary differential equation

$$\hat{v}_{tt} + |\xi|^2 \hat{v} + \tfrac{\mu}{1+t}\hat{v}_t = 0. \tag{2.2}$$

Following [TT80], we use the relation of this differential equation to Bessel's equation to construct a system of linearly independent solutions.

We substitute $\tau = (1+t)|\xi|$ and obtain

$$\frac{\mathrm{d}}{\mathrm{d}t} = \frac{\mathrm{d}\tau}{\mathrm{d}t}\frac{\mathrm{d}}{\mathrm{d}\tau} = |\xi|\frac{\mathrm{d}}{\mathrm{d}\tau} \tag{2.3a}$$

$$\frac{\mathrm{d}^2}{\mathrm{d}t^2} = |\xi|^2 \frac{\mathrm{d}^2}{\mathrm{d}\tau^2} \tag{2.3b}$$

and hence the differential equation

$$|\xi|^2 \hat{v}_{\tau\tau} + |\xi|^2 \hat{v} + \tfrac{\mu}{1+t}|\xi|\hat{v}_\tau = 0.$$

After multiplication with $|\xi|^{-2}$ we get the $|\xi|$-independent differential equation

$$\hat{v}_{\tau\tau} + \tfrac{\mu}{\tau}\hat{v}_\tau + \hat{v} = 0. \tag{2.4}$$

If we make the ansatz $\hat{v} = \tau^\rho w(\tau)$ this leads to

$$\begin{aligned}0 =& \rho(\rho-1)\tau^{\rho-2}w + 2\rho\tau^{\rho-1}w' + \tau^\rho w'' \\&+ \frac{\mu}{\tau}(\rho\tau^{\rho-1}w + \tau^\rho w') + \tau^\rho w \\&+\tau^{\rho-2}\big(\tau^2 w'' + (\mu+2\rho)\tau w' + (\tau^2 + \rho(\rho-1+\mu))w\big),\end{aligned}$$

i.e. by the choice of $\mu + 2\rho = 1$,

$$\rho = -\frac{\mu-1}{2}, \tag{2.5}$$

and hence $\rho - 1 + \mu = -\rho$ we get Bessel's differential equation

$$\tau^2 w'' + \tau w' + (\tau^2 - \rho^2)w = 0 \tag{2.6}$$

of order $\pm \rho$. A system of linearly independent solutions of (2.6) is given by the pair of Hankel functions $\mathcal{H}_\rho^\pm(\tau)$. For details we refer to the treatment in [Wat22].

Hence

$$w_+(\tau) = \tau^\rho \mathcal{H}_\rho^+(\tau), \qquad w_-(\tau) = \tau^\rho \mathcal{H}_\rho^-(\tau), \tag{2.7}$$

with $\rho$ determined by (2.5) gives a pair of linearly independent solutions of (2.4) and after performing the substitution $\tau = (1+t)|\xi|$ of (2.2).

**Representation of the Fourier multiplier.** We are interested in a pair of special solutions $\Phi_1(t, t_0, \xi)$, $\Phi_2(t, t_0, \xi)$ of (2.2) to initial conditions

$$\Phi_1(t_0, t_0, \xi) = 1, \qquad \partial_t \Phi_1(t_0, t_0, \xi) = 0, \tag{2.8a}$$
$$\Phi_2(t_0, t_0, \xi) = 0, \qquad \partial_t \Phi_2(t_0, t_0, \xi) = 1, \tag{2.8b}$$

where the parameter $t_0 > -1$ describes the time level, where the initial data are given. We collect these $\Phi_i$ in the matrix

$$\Phi(t, t_0, \xi) = \begin{pmatrix} \Phi_1(t, t_0, \xi) & \Phi_2(t, t_0, \xi) \\ \partial_t \Phi_1(t, t_0, \xi) & \partial_t \Phi_2(t, t_0, \xi) \end{pmatrix}. \tag{2.9}$$

Then for the vector $(\hat{v}(t, \cdot), \hat{v}_t(t, \cdot))^T$ we obtain the relation

$$\begin{pmatrix} \hat{v}(t, \cdot) \\ \hat{v}_t(t, \cdot) \end{pmatrix} = \Phi(t, t_0, \cdot) \begin{pmatrix} \hat{v}(t_0, \cdot) \\ \hat{v}_t(t_0, \cdot) \end{pmatrix}. \tag{2.10}$$

The matrix $\Phi$ is called the fundamental matrix of (2.2).

For $w_\pm(t, \xi)$ we have the following initial values

$$w_+(t_0, \xi) = (1 + t_0)^\rho |\xi|^\rho \mathcal{H}_\rho^+\big((1 + t_0)|\xi|\big), \tag{2.11a}$$
$$\partial_t w_+(t_0, \xi) = (1 + t_0)^\rho |\xi|^{\rho+1} \mathcal{H}_{\rho-1}^+\big((1 + t_0)|\xi|\big), \tag{2.11b}$$
$$w_-(t_0, \xi) = (1 + t_0)^\rho |\xi|^\rho \mathcal{H}_\rho^-\big((1 + t_0)|\xi|\big), \tag{2.11c}$$
$$\partial_t w_-(t_0, \xi) = (1 + t_0)^\rho |\xi|^{\rho+1} \mathcal{H}_{\rho-1}^-\big((1 + t_0)|\xi|\big), \tag{2.11d}$$

which follow straightforward from the recursion formulae for Bessel functions [Wat22, §3.6]. For instance we have

$$\frac{\mathrm{d}}{\mathrm{d}t}\left[(1+t)^\rho|\xi|^\rho\mathcal{H}_\rho^+((1+t)|\xi|)\right]\Big|_{t=t_0}$$
$$= \rho(1+t)^{\rho-1}|\xi|^\rho\mathcal{H}_\rho^+((1+t)|\xi|)\big|_{t=t_0} + (1+t)^\rho|\xi|^\rho(\mathcal{H}_\rho^+)'((1+t)|\xi|)|\xi|\big|_{t=t_0}$$
$$= (1+t_0)^\rho|\xi|^{\rho+1}\left(\frac{\rho}{|\xi|}\mathcal{H}_\rho^+\big((1+t_0)|\xi|\big) + (\mathcal{H}_\rho^+)'\big((1+t_0)|\xi|\big)\right)$$
$$= (1+t_0)^\rho|\xi|^{\rho+1}\mathcal{H}_{\rho-1}^+\big((1+t_0)|\xi|\big).$$

From these initial values we determine the constants $C_{i\pm}(t_0,\xi)$ in

$$\Phi_i(t,t_0,\xi) = C_{i+}(t_0,\xi)w_+(t,\xi) + C_{i-}(t_0,\xi)w_-(t,\xi), \qquad i=1,2, \quad (2.12)$$

such that (2.8) holds. This means we have to satisfy

$$\begin{pmatrix} w_+(t_0,\xi) & w_-(t_0,\xi) \\ \partial_t w_+(t_0,\xi) & \partial_t w_-(t_0,\xi) \end{pmatrix}\begin{pmatrix} C_{1+}(t_0,\xi) & C_{2+}(t_0,\xi) \\ C_{1-}(t_0,\xi) & C_{2-}(t_0,\xi) \end{pmatrix} = I. \quad (2.13)$$

Hence we have

$$\begin{pmatrix} C_{1+}(t_0,\xi) & C_{2+}(t_0,\xi) \\ C_{1-}(t_0,\xi) & C_{2-}(t_0,\xi) \end{pmatrix} = \frac{1}{\det W(t_0)}\begin{pmatrix} \partial_t w_-(t_0,\xi) & -w_-(t_0,\xi) \\ -\partial_t w_+(t_0,\xi) & w_+(t_0,\xi) \end{pmatrix}, \quad (2.14)$$

where

$$\det W(t_0) = \det\begin{pmatrix} w_+(t_0,\xi) & w_-(t_0,\xi) \\ \partial_t w_+(t_0,\xi) & \partial_t w_-(t_0,\xi) \end{pmatrix}$$
$$= (1+t_0)^{2\rho}|\xi|^{2\rho+1}\det\begin{pmatrix} \mathcal{H}_\rho^+\big((1+t_0)|\xi|\big) & (\mathcal{H}_\rho^+)'\big((1+t_0)|\xi|\big) \\ \mathcal{H}_\rho^-\big((1+t_0)|\xi|\big) & (\mathcal{H}_\rho^-)'\big((1+t_0)|\xi|\big) \end{pmatrix}^T$$
$$= -\frac{4i}{\pi}|\xi|^{2\rho}(1+t_0)^{2\rho-1}, \quad (2.15)$$

because the Wronskian $\mathcal{W}$ of the Hankel functions satisfies

$$\mathcal{W}(\mathcal{H}_\rho^+(z),\mathcal{H}_\rho^-(z)) = \det\begin{pmatrix} \mathcal{H}_\rho^+ & (\mathcal{H}_\rho^+)' \\ \mathcal{H}_\rho^- & (\mathcal{H}_\rho^-)' \end{pmatrix} = -\frac{4i}{\pi z}, \quad (2.16)$$

see [Wat22, §3.63].

Hence we obtain for the fundamental solution

$$\begin{aligned}\Phi_1(t,t_0,\xi) &= C_{1+}(t_0,\xi)w_+(t,\xi) + C_{1-}(t_0,\xi)w_-(t,\xi)\\ &= \frac{i\pi}{4}|\xi|^{-2\rho}(1+t_0)^{1-2\rho}\\ &\quad \Big\{(1+t_0)^\rho|\xi|^{\rho+1}\mathcal{H}^-_{\rho-1}((1+t_0)|\xi|)\,(1+t)^\rho|\xi|^\rho\mathcal{H}^+_\rho((1+t)|\xi|)\\ &\quad - (1+t_0)^\rho|\xi|^{\rho+1}\mathcal{H}^+_{\rho-1}((1+t_0)|\xi|)\,(1+t)^\rho|\xi|^\rho\mathcal{H}^-_\rho((1+t)|\xi|)\Big\}\\ &= \frac{i\pi}{4}|\xi|\frac{(1+t)^\rho}{(1+t_0)^{\rho-1}}\Big\{\mathcal{H}^-_{\rho-1}((1+t_0)|\xi|)\mathcal{H}^+_\rho((1+t)|\xi|)\\ &\quad - \mathcal{H}^+_{\rho-1}((1+t_0)|\xi|)\mathcal{H}^-_\rho((1+t)|\xi|)\Big\} \end{aligned} \qquad (2.17a)$$

and similar

$$\begin{aligned}\Phi_2(t,t_0,\xi) &= -\frac{i\pi}{4}\frac{(1+t)^\rho}{(1+t_0)^{\rho-1}}\Big\{\mathcal{H}^-_\rho((1+t_0)|\xi|)\mathcal{H}^+_\rho((1+t)|\xi|)\\ &\quad - \mathcal{H}^+_\rho((1+t_0)|\xi|)\mathcal{H}^-_\rho((1+t)|\xi|)\Big\}. \end{aligned} \qquad (2.17b)$$

We collect the results in the following theorem.

**Theorem 2.1.** *Assume $v = v(t,x)$ solves Cauchy problem (1.1). Then the Fourier transform $\hat{v}(t,\xi)$ satisfies*

$$\hat{v}(t,\xi) = \sum_{j=1,2} \Phi_j(t,0,\xi)\hat{v}_j(\xi),$$

*where the multipliers $\Phi_j$ are given by*

$$\Phi_1(t,t_0,\xi) = \frac{i\pi}{4}|\xi|\frac{(1+t)^\rho}{(1+t_0)^{\rho-1}}\begin{vmatrix}\mathcal{H}^-_{\rho-1}((1+t_0)|\xi|) & \mathcal{H}^-_\rho((1+t)|\xi|)\\ \mathcal{H}^+_{\rho-1}((1+t_0)|\xi|) & \mathcal{H}^+_\rho((1+t)|\xi|)\end{vmatrix}$$

*and*

$$\Phi_2(t,t_0,\xi) = -\frac{i\pi}{4}\frac{(1+t)^\rho}{(1+t_0)^{\rho-1}}\begin{vmatrix}\mathcal{H}^-_\rho((1+t_0)|\xi|) & \mathcal{H}^-_\rho((1+t)|\xi|)\\ \mathcal{H}^+_\rho((1+t_0)|\xi|) & \mathcal{H}^+_\rho((1+t)|\xi|)\end{vmatrix}.$$

**Time derivatives.** Next we need the time derivatives of these functions. Derivation with respect to $t$ leads to

$$\partial_t \Phi_1(t, t_0, \xi) = \frac{i\pi}{4}|\xi|\rho \frac{(1+t)^{\rho-1}}{(1+t_0)^{\rho-1}} \begin{vmatrix} \mathcal{H}_{\rho-1}^-((1+t_0)|\xi|) & \mathcal{H}_\rho^-((1+t)|\xi|) \\ \mathcal{H}_{\rho-1}^+((1+t_0)|\xi|) & \mathcal{H}_\rho^+((1+t)|\xi|) \end{vmatrix}$$

$$+ \frac{i\pi}{4}|\xi|\frac{(1+t)^\rho}{(1+t_0)^{\rho-1}} \begin{vmatrix} \mathcal{H}_{\rho-1}^-((1+t_0)|\xi|) & (\mathcal{H}_\rho^-)'((1+t)|\xi|) \\ \mathcal{H}_{\rho-1}^+((1+t_0)|\xi|) & (\mathcal{H}_\rho^+)'((1+t)|\xi|) \end{vmatrix} |\xi|$$

$$= \frac{i\pi}{4}|\xi|^2 \frac{(1+t)^\rho}{(1+t_0)^{\rho-1}} \begin{vmatrix} \mathcal{H}_{\rho-1}^-((1+t_0)|\xi|) & \mathcal{H}_{\rho-1}^-((1+t)|\xi|) \\ \mathcal{H}_{\rho-1}^+((1+t_0)|\xi|) & \mathcal{H}_{\rho-1}^+((1+t)|\xi|) \end{vmatrix}$$

using $\rho \mathcal{H}_\rho^+(z) + z (\mathcal{H}_\rho^+)'(z) = z \mathcal{H}_{\rho-1}^+(z)$ and similarly for the second multiplier.

**Corollary 2.2.** *The time derivatives of the multipliers $\Phi_j$ are given by*

$$\partial_t \Phi_1(t, t_0, \xi) = \frac{i\pi}{4}|\xi|^2 \frac{(1+t)^\rho}{(1+t_0)^{\rho-1}} \begin{vmatrix} \mathcal{H}_{\rho-1}^-((1+t_0)|\xi|) & \mathcal{H}_{\rho-1}^-((1+t)|\xi|) \\ \mathcal{H}_{\rho-1}^+((1+t_0)|\xi|) & \mathcal{H}_{\rho-1}^+((1+t)|\xi|) \end{vmatrix}$$

*and*

$$\partial_t \Phi_2(t, t_0, \xi) = -\frac{i\pi}{4}|\xi|\frac{(1+t)^\rho}{(1+t_0)^{\rho-1}} \begin{vmatrix} \mathcal{H}_\rho^-((1+t_0)|\xi|) & \mathcal{H}_{\rho-1}^-((1+t)|\xi|) \\ \mathcal{H}_\rho^+((1+t_0)|\xi|) & \mathcal{H}_{\rho-1}^+((1+t)|\xi|) \end{vmatrix}.$$

It is possible to obtain a similar expression for higher order time derivatives in the same way.

**Representation by real-valued functions.** If we use the definition of $\mathcal{H}_\rho^\pm$ by the real-valued Bessel and Weber functions [Wat22, §3.6]

$$\mathcal{H}_\rho^\pm(z) = \mathcal{J}_\rho(z) \pm i \mathcal{Y}_\rho(z), \tag{2.18}$$

we obtain an alternative characterization of $\Phi$ by real-valued functions

$$\Phi_1(t, t_0, \xi) = -\frac{\pi}{2}|\xi|\frac{(1+t)^\rho}{(1+t_0)^{\rho-1}} \begin{vmatrix} \mathcal{J}_{\rho-1}((1+t_0)|\xi|) & \mathcal{J}_\rho((1+t)|\xi|) \\ \mathcal{Y}_{\rho-1}((1+t_0)|\xi|) & \mathcal{Y}_\rho((1+t)|\xi|) \end{vmatrix}$$

(2.19a)

$$\Phi_2(t, t_0, \xi) = \frac{\pi}{2}\frac{(1+t)^\rho}{(1+t_0)^{\rho-1}} \begin{vmatrix} \mathcal{J}_\rho((1+t_0)|\xi|) & \mathcal{J}_\rho((1+t)|\xi|) \\ \mathcal{Y}_\rho((1+t_0)|\xi|) & \mathcal{Y}_\rho((1+t)|\xi|) \end{vmatrix} \tag{2.19b}$$

$$\partial_t \Phi_1(t, t_0, \xi) = -\frac{\pi}{2}|\xi|^2 \frac{(1+t)^\rho}{(1+t_0)^{\rho-1}} \begin{vmatrix} \mathcal{J}_{\rho-1}((1+t_0)|\xi|) & \mathcal{J}_{\rho-1}((1+t)|\xi|) \\ \mathcal{Y}_{\rho-1}((1+t_0)|\xi|) & \mathcal{Y}_{\rho-1}((1+t)|\xi|) \end{vmatrix}$$

(2.19c)

$$\partial_t \Phi_2(t, t_0, \xi) = \frac{\pi}{2}|\xi|\frac{(1+t)^\rho}{(1+t_0)^{\rho-1}} \begin{vmatrix} \mathcal{J}_\rho((1+t_0)|\xi|) & \mathcal{J}_{\rho-1}((1+t)|\xi|) \\ \mathcal{Y}_\rho((1+t_0)|\xi|) & \mathcal{Y}_{\rho-1}((1+t)|\xi|) \end{vmatrix}.$$

(2.19d)

The representation simplifies in the case of non-integral $\rho$ using

$$\mathcal{Y}_\rho(z) = \cot(\rho\pi)\mathcal{J}_\rho(z) - \csc(\rho\pi)\mathcal{J}_{-\rho}(z) \tag{2.20}$$

and the periodicity properties of the trigonometric functions to

$$\Phi_1(t, t_0, \xi) = -\frac{\pi}{2}\csc(\rho\pi)|\xi|\frac{(1+t)^\rho}{(1+t_0)^{\rho-1}}$$
$$\begin{vmatrix} \mathcal{J}_{-\rho+1}((1+t_0)|\xi|) & \mathcal{J}_{-\rho}((1+t)|\xi|) \\ -\mathcal{J}_{\rho-1}((1+t_0)|\xi|) & \mathcal{J}_\rho((1+t)|\xi|) \end{vmatrix} \tag{2.21a}$$

$$\Phi_2(t, t_0, \xi) = \frac{\pi}{2}\csc(\rho\pi)\frac{(1+t)^\rho}{(1+t_0)^{\rho-1}}$$
$$\begin{vmatrix} \mathcal{J}_{-\rho}((1+t_0)|\xi|) & \mathcal{J}_{-\rho}((1+t)|\xi|) \\ \mathcal{J}_\rho((1+t_0)|\xi|) & \mathcal{J}_\rho((1+t)|\xi|) \end{vmatrix} \tag{2.21b}$$

$$\partial_t \Phi_1(t, t_0, \xi) = -\frac{\pi}{2}\csc(\rho\pi)|\xi|^2\frac{(1+t)^\rho}{(1+t_0)^{\rho-1}}$$
$$\begin{vmatrix} \mathcal{J}_{-\rho+1}((1+t_0)|\xi|) & \mathcal{J}_{-\rho+1}((1+t)|\xi|) \\ \mathcal{J}_{\rho-1}((1+t_0)|\xi|) & \mathcal{J}_{\rho-1}((1+t)|\xi|) \end{vmatrix} \tag{2.21c}$$

$$\partial_t \Phi_2(t, t_0, \xi) = \frac{\pi}{2}\csc(\rho\pi)|\xi|\frac{(1+t)^\rho}{(1+t_0)^{\rho-1}}$$
$$\begin{vmatrix} \mathcal{J}_{-\rho}((1+t_0)|\xi|) & \mathcal{J}_{-\rho+1}((1+t)|\xi|) \\ -\mathcal{J}_\rho((1+t_0)|\xi|) & \mathcal{J}_{\rho-1}((1+t)|\xi|) \end{vmatrix}. \tag{2.21d}$$

In the first and in the last formula we used $\csc(\rho\pi - \pi) = -\csc(\rho\pi)$.

## 3 Estimates

We use the isomorphism

$$\langle D \rangle^s : L_{p,s}(\mathbb{R}^n) \to L_p(\mathbb{R}^n),$$

$\langle D \rangle$ the pseudo-differential operator with symbol $\langle \xi \rangle = \sqrt{1+|\xi|^2}$, to characterize the Sobolev spaces of fractional order[1] over $L_p(\mathbb{R}^n)$, $p \in (1, \infty)$.

The representation of the fundamental matrix $\Phi(t, 0, \xi)$ (as well as the knowledge about strictly hyperbolic problems) imply a natural regularity difference for the data of one Sobolev order. Therefore we define the following two operators corresponding to the Cauchy problem (1.1). On the one hand we are interested in the solution itself. Let

$$\mathbb{S}(t) : (v_1, \langle D \rangle^{-1} v_2) \mapsto v(t, \cdot) \tag{3.1}$$

---
[1] See e.g. [AS61] or [AMS63] for details about Bessel potential spaces.

be the *solution operator*. (Obviously) it maps $\mathbb{S}(t) : L_2(\mathbb{R}^n, \mathbb{R}^2) \to L_2(\mathbb{R}^n)$. On the other hand we are interested in energy estimates. We define the *energy operator*

$$\mathbb{E}(t) : (\langle D \rangle v_1, v_2) \mapsto (\partial_t v(t, \cdot), |D|v(t, \cdot)), \tag{3.2}$$

with $\mathbb{E}(t) : L_2(\mathbb{R}^n, \mathbb{R}^2) \to L_2(\mathbb{R}^n, \mathbb{R}^2)$.

We will give norm estimates for both operators from $L_p$ scale to $L_q$ with dual indices $p$ and $q$, i.e. $p + q = pq$.

*Remark* 3.1. The norm of the data "contains" $\langle D \rangle v_1$ and not only the homogeneous component $|D|v_1$. The dissipative term has influences like a (possibly negative) mass term and thus transfers information of the size of $v$ to the energy $v'$, cf. Appendix A.2 for the transformation of the problem to a Klein-Gordon equation.

**Properties of Bessel functions.** To obtain norm estimates for the operator families $\mathbb{S}(t)$ and $\mathbb{E}(t)$ we have to review some of the main properties of Bessel functions for small and large arguments. For details we refer to [Wat22, §3.52,§10.6 and §7.2].

**Proposition 3.1.**  1. *The function*

$$\Lambda_\nu(\tau) = \tau^{-\nu} \mathcal{J}_\nu(\tau)$$

*is entire in $\nu$ and $\tau$. Furthermore, $\Lambda_\nu(0) \neq 0$.*

2. *Weber's function $\mathcal{Y}_n(\tau)$ satisfies for integral $n$*

$$\mathcal{Y}_n(\tau) = \frac{2}{\pi} \mathcal{J}_n(\tau) \log \tau + A_n(\tau),$$

*where $\tau^n A_n(\tau)$ is entire and non-zero for $\tau = 0$.*

3. *The Hankel functions $\mathcal{H}_\nu^\pm(\tau)$ with $\tau \geq K$ can be written as*

$$\mathcal{H}_\nu^\pm(\tau) = e^{\pm i\tau} a_\nu^\pm(\tau),$$

*where $a_\nu^\pm \in S^{-\frac{1}{2}}(K, \infty)$ is a classical symbol of order $-1/2$.*

4. *For small arguments, $0 < \tau \leq c < 1$, we have*

$$\left|\mathcal{H}_\nu^\pm(\tau)\right| \lesssim \begin{cases} \tau^{-|\nu|} & , \nu \neq 0, \\ -\log \tau & , \nu = 0. \end{cases}$$

## 3.1  $L_2$ and $L_p$–$L_q$ Estimates for a Model Operator

**The model operator.**  Due to the special structure of the fundamental matrix $\Phi(t,t_0,\xi)$ and hence of the multipliers of $\mathbb{S}(t)$ and $\mathbb{E}(t)$, we consider the time dependent model multiplier

$$\Psi_{k,s,\rho,\delta}(t,\xi) = |\xi|^k \langle\xi\rangle^{s+1-k} \begin{vmatrix} \mathcal{H}_\rho^-(|\xi|) & \mathcal{H}_{\rho+\delta}^-((1+t)|\xi|) \\ \mathcal{H}_\rho^+(|\xi|) & \mathcal{H}_{\rho+\delta}^+((1+t)|\xi|) \end{vmatrix} \qquad (3.3)$$

parameterized by $k, s, \rho, \delta \in \mathbb{R}$.

Again we can write $\Psi_{k,s,\rho,\delta}(t,|\xi|)$ in terms of the real-valued Bessel functions of first and second kind. Similar to (2.19) and (2.21) we have

$$\Psi_{k,s,\rho,\delta}(t,\xi) = 2i|\xi|^k \langle\xi\rangle^{s+1-k} \begin{vmatrix} \mathcal{J}_\rho(|\xi|) & \mathcal{J}_{\rho+\delta}((1+t)|\xi|) \\ \mathcal{Y}_\rho(|\xi|) & \mathcal{Y}_{\rho+\delta}((1+t)|\xi|) \end{vmatrix} \qquad (3.4a)$$

$$= 2i \csc(\rho\pi) |\xi|^k \langle\xi\rangle^{s+1-k} \begin{vmatrix} \mathcal{J}_{-\rho}(|\xi|) & \mathcal{J}_{-\rho-\delta}((1+t)|\xi|) \\ (-1)^\delta \mathcal{J}_\rho(|\xi|) & \mathcal{J}_{\rho+\delta}((1+t)|\xi|) \end{vmatrix}, \qquad (3.4b)$$

the last line holds for $\rho \notin \mathbb{Z}$ and $\rho+\delta \notin \mathbb{Z}$.

**$L_2$-estimates.**  The structure of the multiplier hints to a decomposition of the phase space $\mathbb{R}_+ \times \mathbb{R}^n_\xi$ into three regions

$$Z_1 := \{|\xi| \geq K\} \qquad (3.5a)$$
$$Z_2 := \{|\xi| \leq K \leq (1+t)|\xi|\} \qquad (3.5b)$$
$$Z_3 := \{(1+t)|\xi| \leq K\} \qquad (3.5c)$$

with a fixed parameter $K > 0$.

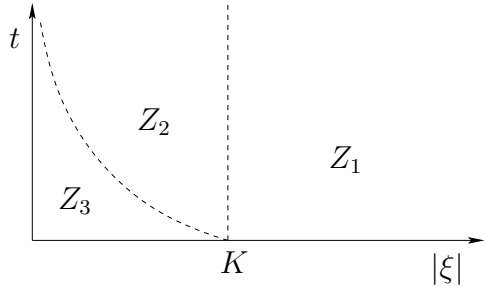

We prove

$$\Psi_{k,0,\rho,\delta}(t,\cdot) \in L_\infty(\mathbb{R}^n) \qquad (3.6)$$

under certain conditions on $k, \rho$ and $\delta$ and derive an estimate for

$$\sup_{\xi \in \mathbb{R}} \big| \langle \xi \rangle^{-s} \Psi_{k,s,\rho,\delta}(t,\xi) \big|. \tag{3.7}$$

From the isometry $M_2^2 \simeq L_\infty$ (see e.g. [Hör60, Theorem 1.5], in general $M_p^q$ denotes the space of all Fourier multipliers mapping $L_p$ into $L_q$) this gives a norm estimate for the operator corresponding to the multiplier (3.3) as operator $H^s \to L_2$.

**Lemma 3.2.** *It holds $\Psi_{k,s,\rho,\delta}(t,\cdot) \in L_\infty(\mathbb{R}^n)$ for all $t$ if and only if $s \leq 0$ and $k \geq |\delta|$. Furthermore the estimate*

$$\big\| \Psi_{k,s,\rho,\delta}(t,\cdot) \big\|_\infty \lesssim \begin{cases} (1+t)^{-\frac{1}{2}} & , \rho \neq 0, |\rho| - k \leq -\frac{1}{2}, \\ (1+t)^{|\rho|-k} & , \rho \neq 0, |\rho| - k \geq -\frac{1}{2}, \\ (1+t)^{-k} \log(e+t) & , \rho = 0, k \leq \frac{1}{2}. \end{cases}$$

*is valid.*

*Proof.* We subdivide the proof into three parts corresponding to the three zones $Z_1$, $Z_2$ and $Z_3$.

$\boxed{Z_1}$ We use Proposition 3.1.3 together with the definition of the zone $K \leq |\xi| \leq (1+t)|\xi|$. Thus the multiplier is bounded in $Z_1$ iff $s \leq 0$. It satisfies

$$\big| \Psi_{k,s,\rho,\delta}(t,\xi) \big| \lesssim (1+t)^{-\frac{1}{2}}$$

under this assumption.

$\boxed{Z_2}$ For $\rho \neq 0$ we can use 3.1.4 to conclude

$$\big| \Psi_{k,s,\rho,\delta}(t,\xi) \big| \lesssim |\xi|^{k-|\rho|}(1+t)^{-\frac{1}{2}}|\xi|^{-\frac{1}{2}}$$

$$\lesssim \begin{cases} (1+t)^{-\frac{1}{2}} & , |\rho| - k \leq -\frac{1}{2}, \\ (1+t)^{|\rho|-k} & , |\rho| - k \geq -\frac{1}{2}. \end{cases}$$

For $\rho = 0$ we have to modify this estimate by the log term

$$\big| \Psi_{k,s,0,\delta}(t,\xi) \big| \lesssim |\xi|^k \log \frac{2K}{|\xi|} (1+t)^{-\frac{1}{2}}|\xi|^{-\frac{1}{2}}$$

$$\lesssim \begin{cases} (1+t)^{-\frac{1}{2}} & , k > \frac{1}{2}, \\ (1+t)^{-k} \log(e+t) & , k \leq \frac{1}{2}. \end{cases}$$

$\boxed{Z_3}$ In this zone we use the representation of $\Psi_{k,s,\rho,\delta}(t,\xi)$ in terms of real-valued functions given by (3.4). For non-integral $\rho$ and $\rho+\delta$ we can use the representation by Bessel functions of first kind to conclude[2]

$$\begin{aligned}
\big|\Psi_{k,s,\rho,\delta}(t,\xi)\big| & \\
&\lesssim \left||\xi|^\rho \mathcal{J}_{-\rho}(|\xi|)\,\big((1+t)|\xi|\big)^{-\rho-\delta}\mathcal{J}_{\rho+\delta}\big((1+t)|\xi|\big)\,(1+t)^{\rho+\delta}|\xi|^{k+\delta}\right| \\
&\quad + \left||\xi|^{-\rho}\mathcal{J}_\rho(|\xi|)\,\big((1+t)|\xi|\big)^{\rho+\delta}\mathcal{J}_{-\rho-\delta}\big((1+t)|\xi|\big)\,(1+t)^{-\rho-\delta}|\xi|^{k-\delta}\right| \\
&\lesssim (1+t)^{|\rho|-k}.
\end{aligned}$$

The condition $k \geq |\delta|$ is necessary and sufficient for the boundedness in $\xi$.

For integral values of $\rho$ or $\rho+\delta$ we use Weber's functions and Proposition 3.1.2. We sketch the estimate if both $\rho$ and $\rho+\delta$ are integral. Then we have

$$\begin{aligned}
\Psi_{k,s,\rho,\delta} &= -\frac{4i}{\pi}|\xi|^k \langle\xi\rangle^{s+1-k} \log(1+t)\mathcal{J}_\rho(|\xi|)\mathcal{J}_{\rho+\delta}\big((1+t)|\xi|\big) \\
&\quad + 2i|\xi|^k \langle\xi\rangle^{s+1-k} \begin{vmatrix} \mathcal{J}_\rho(|\xi|) & \mathcal{J}_{\rho+\delta}\big((1+t)|\xi|\big) \\ A_\rho(|\xi|) & A_{\rho+\delta}\big((1+t)|\xi|\big) \end{vmatrix}
\end{aligned}$$

and hence

$$\begin{aligned}
\big|\Psi_{k,s,\rho,\delta}(t,\xi)\big| &\lesssim \log(e+t)\,(1+t)^{-|\rho|-k} + (1+t)^{|\rho|-k} \\
&\lesssim \begin{cases} (1+t)^{|\rho|-k} & ,\rho \neq 0, \\ (1+t)^{-k}\log(e+t) & ,\rho = 0. \end{cases}
\end{aligned}$$

If only one of both indices is integral we have to mix the representations.[3] $\square$

$L_p$–$L_q$-**estimates.** We consider estimates for the model operator

$$u(x) \mapsto \mathcal{F}^{-1}\big[\Psi_{k,s,\rho,\delta}(t,\xi)\hat{u}(\xi)\big](x) \qquad (3.8)$$

from $L_{p,r}(\mathbb{R}^n)$ to $L_q(\mathbb{R}^n)$, $(p,q)$ a dual pair and $r$ a sufficiently high regularity.

---

[2] By estimating the difference structure of the multiplier by triangle inequality we do not loose information. For non-integral $\rho$ the leading terms of the series expansions do not cancel.

[3] For later reference we use only integral values of $\delta$.

**Lemma 3.3.** *Assume $p \in (1, 2]$, $q$ such that $pq = p + q$. Let further $k \geq |\delta|$. Then the model operator (3.8) satisfies the norm estimate*

$$||\Psi_{k,s,\rho,\delta}(t, \mathrm{D})||_{L_{p,r} \to L_q}$$
$$\lesssim \begin{cases} (1+t)^{-\frac{n-1}{2}\left(\frac{1}{p}-\frac{1}{q}\right)-\frac{1}{2}} & , d > \frac{1}{2} \\ (1+t)^{-n\left(\frac{1}{p}-\frac{1}{q}\right)+|\rho|-k} & , \rho \neq 0, d \leq \frac{1}{2} \\ (1+t)^{-n\left(\frac{1}{p}-\frac{1}{q}\right)+\theta\epsilon-k}(\log(e+t))^{1-\theta} & , \rho = 0, d < \frac{1}{2} + \epsilon, \epsilon > 0 \end{cases}$$

*for $d = \frac{n+1}{2}\left(\frac{1}{p} - \frac{1}{q}\right) + k - |\rho|$ and $r = n\left(\frac{1}{p} - \frac{1}{q}\right) + s$. The interpolating constant $\theta$ in the last case is given by $\theta = \frac{n+1}{2\epsilon+1-2k}\left(\frac{1}{p} - \frac{1}{q}\right)$.*

*Proof.* Again we decompose the phase space into three zones. For this we use a smooth cut-off function $\psi \in C^\infty(\mathbb{R}_+)$ with $\psi' \leq 0$, $\psi(r) = 1$ for $r < 1/2$ and $\psi(r) = 0$ for $r > 2$. Using this function we define

$$\phi_1(t, \xi) = 1 - \psi(|\xi|/K),$$
$$\phi_2(t, \xi) = \psi(|\xi|/K)\big(1 - \psi((1+t)|\xi|/K)\big),$$
$$\phi_3(t, \xi) = \psi(|\xi|/K)\psi((1+t)|\xi|/K),$$

such that $\phi_1(t, \xi) + \phi_2(t, \xi) + \phi_3(t, \xi) = 1$ and $\operatorname{supp} \phi_i \sim Z_i$. Thus we can decompose the multiplier $\Psi_{k,s,\rho,\delta}$ into the sum $\sum_{i=1,2,3} \phi_i(t, \xi)\Psi_{k,s,\rho,\delta}(t, \xi)$ and estimate each of the summands. We prove the estimate for $r = 0$, i.e. we restrict the proof to the corresponding value

$$s = -n\left(\frac{1}{p} - \frac{1}{q}\right).$$

$\boxed{Z_1}$ From Proposition 3.1 we conclude the representation of

$$\phi_1(t, \xi)\Psi_{k,s,\rho,\delta}(t, \xi)$$

as sum of two multipliers of the form

$$e^{\pm t|\xi|} a(|\xi|) b\big((1+t)|\xi|\big)$$

with symbols $a \in S^{s+\frac{1}{2}}(K, \infty)$ and $b \in S^{-\frac{1}{2}}(K, \infty)$. We follow [Bre75] to estimate the corresponding Fourier integral operator. The key tool is a dyadic decomposition together with Littman's Lemma, [Bre75, Lemma 4].

Let $\chi \in C_0^\infty(\mathbb{R}_+)$ be nonnegative with support contained in $[1/2, 2]$ and

$$\sum_{j=-\infty}^{\infty} \chi(2^j r) = 1, \qquad \text{for } r \neq 0. \tag{3.9}$$

Such a function exists from [Hör60, Lemma 2.3]. Let further $\chi_j(\xi) = \chi(2^{-j}\xi/K)$.

We obtain an $L_p$–$L_q$ estimate for this operator by interpolating $L_1$–$L_\infty$ and $L_2$–$L_2$ estimates using [Bre75, Lemma 3]. For this we define

$$I_j = \left|\left| \mathcal{F}^{-1}\left[ \chi_j(\xi) e^{\pm it|\xi|} a(|\xi|) b\big((1+t)|\xi|\big) \right] \right|\right|_\infty \tag{3.10a}$$

and

$$\tilde{I}_j = \left|\left| \chi_j(\xi) e^{\pm it|\xi|} a(|\xi|) b\big((1+t)|\xi|\big) \right|\right|_\infty \tag{3.10b}$$

and estimate these norms. For all $j < 0$ we have $I_j = \tilde{I}_j = 0$.

For $I_j$ we perform the substitution $\xi = 2^j K \eta$ and obtain

$$I_j \leq C 2^{jn} \left|\left| \mathcal{F}^{-1}\left[ e^{\pm it 2^j K \eta} a(2^j K \eta) b((1+t) 2^j K \eta) \chi(|\eta|) \right] \right|\right|_\infty$$

$$\leq C 2^{jn} (1 + 2^j K t)^{-\frac{n-1}{2}} \sum_{|\alpha| \leq M} ||D^\alpha a(2^j K \eta) b((1+t) 2^j K \eta) \chi(|\eta|)||_\infty$$

$$\leq C 2^{jn} (1 + 2^j K t)^{-\frac{n-1}{2}}$$

$$\sum_{|\alpha+\beta| \leq M} \sup_{1/2 \leq |\eta| \leq 2} (2^j K |\eta|)^{s+\frac{1}{2}-|\alpha|} 2^{j|\alpha|} ((1+t) 2^j K |\eta|)^{-\frac{1}{2}-|\beta|} ((1+t) 2^j)^{|\beta|}$$

$$\leq C 2^{j(n+s)} (1 + 2^j K t)^{-\frac{n-1}{2}} (1+t)^{-\frac{1}{2}}$$

by Littmann's Lemma. From $C_K(1+t) \leq (1 + 2^j K t) \leq C'_K 2^j (1+t)$ we get finally

$$I_j \leq C 2^{j(n+s)} (1+t)^{-\frac{n}{2}}. \tag{3.11}$$

For $\tilde{I}_j$ we obtain

$$\tilde{I}_j \leq C \sup_{\eta \in \text{supp } \chi} \phi_1(2^j K \eta) |2^j K \eta|^{s+\frac{1}{2}} |(1+t) 2^j K \eta|^{-\frac{1}{2}}$$

$$\leq C 2^{js} (1+t)^{-\frac{1}{2}}. \tag{3.12}$$

The estimates (3.11) and (3.12) correspond to $L_1$–$L_\infty$ and $L_2$–$L_2$ estimates for the dyadic component of the model operator (3.8). Interpolation

leads to

$$\left\|\mathcal{F}^{-1}\left[\phi_1(t,\xi)\chi_j(\xi)\Psi_{k,s,\rho,\delta}(t,\xi)\hat{u}(\xi)\right]\right\|_q$$
$$\leq C 2^{j\left(n\left(\frac{1}{p}-\frac{1}{q}\right)+s\right)}(1+t)^{-\frac{n-1}{2}\left(\frac{1}{p}-\frac{1}{q}\right)-\frac{1}{2}}\|u\|_p \quad (3.13)$$

for all $p \in (1,2]$, $p+q = pq$. Finally, we use [Bre75, Lemma 2] to conclude for $n\left(\frac{1}{p}-\frac{1}{q}\right)+s \leq 0$ the estimate

$$\left\|\mathcal{F}^{-1}\left[\phi_1(t,\xi)\Psi_{k,s,\rho,\delta}(t,\xi)\hat{u}(\xi)\right]\right\|_q \leq C(1+t)^{-\frac{n-1}{2}\left(\frac{1}{p}-\frac{1}{q}\right)-\frac{1}{2}}\|u\|_p. \quad (3.14)$$

$\boxed{Z_2}$ In this zone we subdivide each summand of the multiplier in two factors

$$|\xi|^{|\rho|+\varepsilon}\mathcal{H}_\rho^\mp(|\xi|)\phi_{21}(\xi), \quad (3.15\text{a})$$
$$\left((1+t)|\xi|\right)^{k-|\rho|-\varepsilon}\mathcal{H}_{\rho+\delta}^\pm\left((1+t)|\xi|\right)\phi_{22}(t,\xi), \quad (3.15\text{b})$$

and the remaining constant $(1+t)^{|\rho|-k+\varepsilon}$, where

$$\phi_{21}(\xi) = \psi(|\xi|/K) \quad \text{and} \quad \phi_{22}(t,\xi) = 1 - \psi((1+t)|\xi|/K) \quad (3.16)$$

such that $\phi_2(t,\xi) = \phi_{21}(\xi)\phi_{22}(t,\xi)$.

The first multiplier is time independent and satisfies

$$|\xi|^{|\rho|+\varepsilon}\mathcal{H}_\rho^\mp(|\xi|)\phi_{21}(\xi) \sim |\xi|^{|\rho|+\varepsilon}\mathcal{H}_{|\rho|}^\mp(|\xi|)\phi_{21}(\xi)$$
$$= (1 \pm i\cot|\rho|\pi)|\xi|^{2|\rho|+\epsilon}|\xi|^{-|\rho|}\mathcal{J}_{|\rho|}(|\xi|)\phi_{21}(\xi)$$
$$\mp i\csc|\rho|\pi\ |\xi|^\epsilon|\xi|^{|\rho|}\mathcal{J}_{-|\rho|}(|\xi|)\phi_{21}(\xi), \qquad \rho \notin \mathbb{Z},$$
$$= |\xi|^{2|\rho|+\varepsilon}|\xi|^{-|\rho|}\mathcal{J}_{|\rho|}(|\xi|)\phi_{21}(\xi)$$
$$\pm i|\xi|^{2|\rho|+\varepsilon}\log|\xi|\ |\xi|^{-|\rho|}\mathcal{J}_{|\rho|}(|\xi|)\phi_{21}(\xi)$$
$$\pm i|\xi|^\varepsilon|\xi|^{|\rho|}A_{|\rho|}(|\xi|)\phi_{21}(\xi), \qquad \rho \in \mathbb{Z},$$

where $\tau^{\mp|\rho|}\mathcal{J}_{\pm\rho}(\tau)$ and $\tau^{|\rho|}A_{|\rho|}(\tau)$ are entire. By $\sim$ we denote equality up to a multiplicative constant here.

From [Ste70, §3 Theorem 3] it follows

$$|\xi|^\varepsilon \phi_{21}(|\xi|) \in M_p^p \qquad \forall \epsilon \geq 0 \quad (3.17\text{a})$$
$$|\xi|^\varepsilon \log|\xi|\ \phi_{21}(|\xi|) \in M_p^p \qquad \forall \epsilon > 0 \quad (3.17\text{b})$$

for all $p \in (1, \infty)$. Thus we conclude with the algebra property of multiplier spaces that the first multiplier belongs to $M_p^p$ for $p \in (1, \infty)$ if $\varepsilon \geq 0$ and $\rho \neq 0$ (or for $\rho = 0$ if $\varepsilon > 0$).

For the second multiplier we prove an $L_p$–$L_q$ estimate. For this we use again a dyadic decompositon. Let $\chi$ be like in the discussion of $Z_1$ and

$$\chi_j(t,\xi) = \chi\big(2^{-j}(1+t)|\xi|/K\big), \tag{3.18a}$$

$$\chi_0(t,\xi) = 1 - \sum_{j>0} \chi_j(t,\xi). \tag{3.18b}$$

We estimate

$$I_j = \left\| \mathcal{F}^{-1}\left[\chi_j(t,\xi)\big((1+t)|\xi|\big)^{k-|\rho|-\varepsilon} \mathcal{H}_{\rho+\delta}^{\pm}\big((1+t)|\xi|\big)\phi_{22}(t,\xi)\right]\right\|_\infty \tag{3.19}$$

and

$$\tilde{I}_j = \left\|\chi_j(t,\xi)\big((1+t)|\xi|\big)^{k-|\rho|-\varepsilon} \mathcal{H}_{\rho+\delta}^{\pm}\big((1+t)|\xi|\big)\phi_{22}(t,\xi)\right\|_\infty. \tag{3.20}$$

For $j > 0$ we have $\phi_{22}(t,\xi) = 1$ on $\mathrm{supp}\,\chi_j$. Hence using the substitution $(1+t)\xi = 2^j K\eta$ we get the estimate

$$\begin{aligned}
I_j &= \left\|\mathcal{F}^{-1}\left[e^{i(1+t)|\xi|} a\big((1+t)\xi\big)\psi_j(t,\xi)\right]\right\|_\infty \\
&\leq C 2^{jn}(1+t)^{-n}\left\|e^{i2^j K\eta} a(2^j K\eta)\psi(|\eta|)\right\|_\infty \\
&\leq C 2^{jn}(1+t)^{-n}(1+2^j K)^{-\frac{n-1}{2}}(2^j K)^{-\frac{1}{2}+k-|\rho|-\varepsilon} \\
&\leq C 2^{j\left(\frac{n+1}{2}+k-|\rho|-\varepsilon-\frac{1}{2}\right)}(1+t)^{-n},
\end{aligned}$$

where $a \in S^{-\frac{1}{2}+k-|\rho|-\varepsilon}$. For $I_0$ we obtain a similar estimate in the same way.

For $\tilde{I}_j$ we have

$$\begin{aligned}
\tilde{I}_j &\leq C \sup_{\eta \in \mathrm{supp}\,\psi} (2^j K\eta)^{k-|\rho|-\varepsilon-\frac{1}{2}} \\
&\leq C 2^{j(k-|\rho|-\varepsilon-\frac{1}{2})}. \tag{3.21}
\end{aligned}$$

Interpolation leads to

$$\left\|\mathcal{F}^{-1}\left[\psi_j(t,\xi)\big((1+t)|\xi|\big)^{k-|\rho|-\varepsilon}\mathcal{H}_{\rho+\delta}^{\pm}\big((1+t)|\xi|\big)\phi_{22}(t,\xi)\hat{u}(\xi)\right]\right\|_q$$
$$\leq C 2^{j\left(\frac{n+1}{2}\left(\frac{1}{p}-\frac{1}{q}\right)+k-|\rho|-\varepsilon-\frac{1}{2}\right)}(1+t)^{-n\left(\frac{1}{p}-\frac{1}{q}\right)}\|u\|_p \tag{3.22}$$

which gives for
$$\varepsilon \geq \frac{n+1}{2}\left(\frac{1}{p}-\frac{1}{q}\right)-\frac{1}{2}+k-|\rho|$$
the $L_p$–$L_q$ estimate
$$\left\|\mathcal{F}^{-1}\left[((1+t)|\xi|)^{k-|\rho|-\varepsilon}\mathcal{H}^{\pm}_{\rho+\delta}((1+t)|\xi|)\phi_{22}(t,\xi)\hat{u}(\xi)\right]\right\|_q$$
$$\leq C(1+t)^{-n\left(\frac{1}{p}-\frac{1}{q}\right)}\|u\|_p. \quad (3.23)$$

The 'regularity' $\varepsilon$ is determined from both multipliers, hence the optimal choice is
$$\varepsilon = \max\left\{0, \tfrac{n+1}{2}\left(\tfrac{1}{p}-\tfrac{1}{q}\right) - \tfrac{1}{2}+k-|\rho|\right\}. \quad (3.24)$$
under the assumption $\rho \neq 0$. For $\rho = 0$ the choice $\epsilon = 0$ is not possible. Therefore for $k \leq \frac{1}{2}$ we have to exclude that case. We postpone this exceptional case.

Multiplication of the multipliers corresponds to a concatenation of the corresponding operators. Hence we have
$$\left\|\mathcal{F}^{-1}\left[\phi_2(t,\xi)\Psi_{k,s,\rho,\delta}(t,\xi)\hat{u}(\xi)\right]\right\|_q$$
$$\leq C(1+t)^{\max\left\{-\frac{n-1}{2}\left(\frac{1}{p}-\frac{1}{q}\right)-\frac{1}{2},\, -n\left(\frac{1}{p}-\frac{1}{q}\right)+|\rho|-k\right\}}\|u\|_p \quad (3.25)$$
for $\rho \neq 0$ or $\rho = 0$ and $k > 1/2$.

$\boxed{Z_2 \text{ for } \rho = 0 \text{ and } k \leq 1/2}$ In this exceptional case we get an estimate for all dual $p$ and $q$ by interpolation. From $\frac{n}{2}+k > 0$ we can follow the previously used reasoning to conclude an estimate for dual $p$, $q$ with
$$\frac{n+1}{2}\left(\frac{1}{p}-\frac{1}{q}\right)-\frac{1}{2}+k > 0.$$

If we interpolate the corresponding result with the previously proven $L_2$ estimate in that zone we obtain
$$\left\|\mathcal{F}^{-1}\left[\phi_2(t,\xi)\Psi_{k,s,0,\delta}(t,\xi)\hat{u}(\xi)\right]\right\|_q \leq C(1+t)^{-n\left(\frac{1}{p}-\frac{1}{q}\right)+\theta\epsilon-k}(\log(e+t))^{1-\theta}\|u\|_p \quad (3.26)$$
for $\frac{n+1}{2}\left(\frac{1}{p}-\frac{1}{q}\right)-\frac{1}{2}+k \leq \epsilon$. The interpolating constant $\theta$ is given by
$$\theta = \frac{n+1}{2\epsilon+1-2k}\left(\frac{1}{p}-\frac{1}{q}\right).$$

$\boxed{Z_3}$ We use the estimate

$$\left|\phi_3(t,\xi)\Psi_{k,s,\rho,\delta}(t,\xi)\right| \lesssim \begin{cases} (1+t)^{|\rho|-k} &, \rho \neq 0, \\ (1+t)^{-k}\log(e+t) &, \rho = 0 \end{cases}$$

together with the definition of the zone $Z_3$ to conclude the estimate

$$\begin{aligned}
||\mathcal{F}[\phi_3(t,\xi)&\Psi_{k,s,\rho,\delta}(t,\xi)\hat{u}(\xi)]||_q \\
&\leq ||\phi_3(t,\cdot)\Psi_{k,s,\rho,\delta}(t,\cdot)\hat{u}(\xi)||_p \\
&\leq ||\phi_3(t,\cdot)||_{1/\left(\frac{1}{p}-\frac{1}{q}\right)}||\hat{u}||_q||\Psi_{k,s,\rho,\delta}(t,\cdot)||_\infty \\
&\lesssim ||u||_p(1+t)^{-n\left(\frac{1}{p}-\frac{1}{q}\right)} \begin{cases} (1+t)^{|\rho|-k} &, \rho \neq 0, \\ (1+t)^{-k}\log(e+t) &, \rho = 0. \end{cases}
\end{aligned}$$

Under the assumptions on $p$ and $q$ this estimate is weaker than the estimates in the zones $Z_1$ and $Z_2$. $\square$

*Remark* 3.2. The zone $Z_1$ determines the necessary regularity of the data, large values of $|\xi|$ are contained in this zone only. This is a difference to the reasoning in [Rei97], [RY98] and [Rei00]. There the pseudo-differential zone contains $\{(t,\xi)|t=0\}$, and hence this zone determines the regularity of the data in this case. This is similar to our consideration if we decompose $Z_1$ in small $t$ and large $t$. For large $t$ we can replace (3.11) by

$$I_j \leq C 2^{j(\frac{n+1}{2}+s)}(1+t)^{-\frac{n}{2}},$$

which gives as conclusion an estimate with less regular data. Hence the necessary regularity is determined by the boarder of $Z_1$.

*Remark* 3.3. It is possible to extend the estimates to the case $p=1$ and $q=\infty$ by Sobolev embedding. Therefore one has to use slightly more regularity, $r = n + \epsilon$.

*Remark* 3.4. In Lemma 3.3 we can weaken the assumption $k \geq |\delta|$ in dependence on the given dual pair $(p,q)$ using the Hardy-Littlewood theorem of fractional integration to estimate in $Z_3$. This will not be necessary for our considerations.

## 3.2 $L_2$ and Energy Estimates

If we compare (2.17) with (3.3) we obtain the representation

$$\Phi_1(t, 0, \xi) = \frac{i\pi}{4}(1+t)^\rho \Psi_{1,0,\rho-1,1}(t, \xi), \tag{3.27a}$$

$$\Phi_2(t, 0, \xi) = -\frac{i\pi}{4}(1+t)^\rho \Psi_{0,-1,\rho,0}(t, \xi), \tag{3.27b}$$

$$\partial_t \Phi_1(t, 0, \xi) = \frac{i\pi}{4}(1+t)^\rho \Psi_{2,1,\rho-1,0}(t, \xi), \tag{3.27c}$$

$$\partial_t \Phi_2(t, 0, \xi) = -\frac{i\pi}{4}(1+t)^\rho \Psi_{1,0,\rho,-1}(t, \xi) \tag{3.27d}$$

of the entries of the fundamental matrix by our model multiplier. Thus we can apply the estimates of Lemma 3.2 to get a-priori estimates for the solution $v = v(t, x)$ of (1.1). This gives

$$||v(t, \cdot)||_2 \lesssim ||v_1||_2 + ||v_2||_{H^{-1}} \begin{cases} (1+t)^{2\rho} &, \rho \in (0, \frac{1}{2}), \\ \log(e+t) &, \rho = 0, \\ 1 &, \rho < 0. \end{cases} \tag{3.28a}$$

For the first derivatives we obtain

$$||v_t(t, \cdot)||_2 \leq C_1(1+t)^{\rho-\frac{1}{2}}||v_1||_{H^1} + C_2(1+t)^{\max\{\rho-\frac{1}{2}, -1\}}||v_2||_2, \tag{3.28b}$$

$$||\nabla v(t, \cdot)||_2 \leq C_1(1+t)^{\rho-\frac{1}{2}}||v_1||_{H^1} + C_2(1+t)^{\max\{\rho-\frac{1}{2}, -1\}}||v_2||_2, \tag{3.28c}$$

which reestablish already known results on the energy decay, [Mat77], [Ues79] and [HN01, Example 2.1].

We collect the above estimates in the following theorem.

**Theorem 3.4.** *1. The solution operator $\mathbb{S}(t)$ satisfies the norm estimate*

$$||\mathbb{S}(t)||_{2 \to 2} \lesssim \begin{cases} (1+t)^{1-\mu} &, \mu \in (0, 1), \\ \log(e+t) &, \mu = 1, \\ 1 &, \mu > 1. \end{cases}$$

*2. The energy operator $\mathbb{E}(t)$ satisfies*

$$||\mathbb{E}(t)||_{2 \to 2} \lesssim \begin{cases} (1+t)^{-\frac{\mu}{2}} &, \mu \in (0, 2], \\ (1+t)^{-1} &, \mu > 2. \end{cases}$$

Especially for $\mu = 2$ this energy estimate is sharper then known estimates. For this special case see also Appendix A.2.

*Remark* 3.5. It should be remarked that it is essential to use the $H^1$-Norm on the right hand side of the energy estimate (or the normalization by $\langle D \rangle^{-1}$ in the definition of $\mathbb{E}(t)$, formula (3.2)). Otherwise we get for the usual energy from Lemma 3.2 only the trivial (and in view of this Lemma also sharp!) estimate

$$E(v;t) \lesssim E(v;0).$$

This implies that information about the size of the datum $v_1$ is necessary for precise a-priori estimates of the energy.

*Remark* 3.6. If one is only interested in the decay order, but not the dependence on the data, one can obtain the slightly better result

$$\lim_{t \to \infty} (1+t)^{\min\{\mu,2\}} E(v;t) = 0.$$

This result follows from the considerations in [HN01].

*Remark* 3.7. We see that the decay rates are determined

- for small $\mu$ by the exterior zones $Z_1$ and $Z_2$ and
- for small $\mu$ by the interior zones $Z_2$ and $Z_3$.

This explains the change in the energy estimate with critical value $\mu = 2$.

## 3.3 $L_p$–$L_q$ Estimates

We apply Lemma 3.3 to the operators $\mathbb{S}(t)$ and $\mathbb{E}(t)$. Therefore we use the representation of the fundamental matrix $\Phi$ by the model multiplier $\Psi_{k,s,\rho,\delta}$ given in the previous section.

**Theorem 3.5.** *1. The solution operator $\mathbb{S}(t)$ satisfies the norm estimate*

$$\|\mathbb{S}(t)\|_{L_{p,r} \to L_q}$$
$$\lesssim \begin{cases} (1+t)^{\max\{-\frac{n-1}{2}\left(\frac{1}{p}-\frac{1}{q}\right)-\frac{\mu}{2},\, -n\left(\frac{1}{p}-\frac{1}{q}\right)+1-\mu\}} & , \mu \in (0,1), \\ (1+t)^{-\frac{n-1}{2}\left(\frac{1}{p}-\frac{1}{q}\right)-\frac{1}{2}} & , \mu = 1, \delta > \frac{1}{2} \\ (1+t)^{-n\left(\frac{1}{p}-\frac{1}{q}\right)+\theta\epsilon}(\log(e+t))^{1-\theta} & , \mu = 1, \delta < \frac{1}{2}+\epsilon, \epsilon > 0 \\ (1+t)^{\max\{-\frac{n-1}{2}\left(\frac{1}{p}-\frac{1}{q}\right)-\frac{\mu}{2},\, -n\left(\frac{1}{p}-\frac{1}{q}\right)\}} & , \mu > 1, \end{cases}$$

*for $p \in (1,2]$, $q$ with $pq = p+q$, $\delta = \frac{n+1}{2}\left(\frac{1}{p}-\frac{1}{q}\right)$ and $r = n\left(\frac{1}{p}-\frac{1}{q}\right)$.*

*The interpolating constant $\theta$ in the third case is given by $\theta = \frac{2d}{2\epsilon+1}$.*

2. The energy operator $\mathbb{E}(t)$ satisfies

$$||\mathbb{E}(t)||_{L_{p,r}\to L_q} \lesssim (1+t)^{\max\{-\frac{n-1}{2}\left(\frac{1}{p}-\frac{1}{q}\right)-\frac{\mu}{2},\, -n\left(\frac{1}{p}-\frac{1}{q}\right)-1\}}$$

for $p \in (1, 2]$, $q$ with $pq = p + q$ and $r = n\left(\frac{1}{p} - \frac{1}{q}\right)$.

*Remark* 3.8. The estimates coincide for $\mu = 0$ (i.e. $\rho = 1/2$) with the well-known $L_p$–$L_q$ estimates for the wave equation [Str70]. For $\mu = 2$ (i.e. $\rho = -1/2$) we can reduce the Cauchy problem (1.1) to the Cauchy problem for the wave equation, see in Appendix A.2.

*Remark* 3.9. Finally, these estimates are related to $L_p$–$L_q$ estimates for equations with increasing in time coefficients. By a change of coordinates (already mentioned in the introduction) we can transform problem (1.1) to

$$v_{tt} - \lambda^2(t)\Delta v = 0$$

with $\lambda(t) = (1+t)^\ell$, $\ell > 0$ for $\mu \in (0,1)$ and $\lambda(t) = e^t$ for $\mu = 1$. In [Rei97] and [Gal00] $L_p$–$L_q$ estimates for these problems with increasing in time coefficients were studied. The estimates obtained here for $\mu \in (0,1]$ imply the estimates given in these papers.

## 4 Inhomogeneous Problems

We conclude the discussion of solution representations and a-priori estimates with remarks to the inhomogeneous Cauchy problem

$$\Box v + \frac{\mu}{1+t}v_t = f, \quad v(0,\cdot) = v_1, \quad v_t(0,\cdot) = v_2 \tag{4.1}$$

to data $v_1, v_2 \in \mathcal{S}(\mathbb{R}^n)$ and a right-hand side $f(t,x) \in C^\infty(\mathbb{R}_+, \mathcal{S}(\mathbb{R}^n))$.

### 4.1 Solution Representation, Duhamel's Principle

If we have a nontrivial right-hand side in (4.1) we need entries of the full fundamental matrix $\Phi(t, t_0, \xi)$ to represent the solution.

To obtain the representation we again apply a partial Fourier transform to (4.1). This yields

$$\hat{v}_{tt} + \frac{\mu}{1+t}\hat{v}_t + |\xi|^2 \hat{v} = \hat{f} \tag{4.2}$$

with data $\hat{v}(0,\cdot) = \hat{v}_1$ and $\hat{v}_t(0,\cdot) = \hat{v}_2$. If $f \equiv 0$ we have

$$\hat{v}(t,\xi) = \sum_{i=1,2} \Phi_i(t,0,\xi)\hat{v}_i \tag{4.3}$$

from the discussion of Section 2. Now we assume $\hat{v}_i \equiv 0$. We claim that in this case
$$\hat{v}(t,\xi) = \int_0^t \Phi_2(t,\tau,\xi)\hat{f}(\tau,\xi)\mathrm{d}\tau. \tag{4.4}$$

Indeed we have
$$\partial_t \int_0^t \Phi_2(t,\tau,\xi)\hat{f}(\tau,\xi)\mathrm{d}\tau$$
$$= \Phi_2(t,t,\xi)\hat{f}(t,\xi) + \int_0^t \partial_t\Phi_2(t,\tau,\xi)\hat{f}(\tau,\xi)\mathrm{d}\tau$$
$$= \int_0^t \partial_t\Phi_2(t,\tau,\xi)\hat{f}(\tau,\xi)\mathrm{d}\tau$$

by the initial condition $\Phi_2(t,t,\xi) = 0$ and
$$\partial_t^2 \int_0^t \Phi_2(t,\tau,\xi)\hat{f}(\tau,\xi)\mathrm{d}\tau$$
$$= (\partial_t\Phi_2)(t,t,\xi)\hat{f}(t,\xi) + \int_0^t \partial_t^2\Phi_2(t,\tau,\xi)\hat{f}(\tau,\xi)\mathrm{d}\tau$$
$$= \hat{f}(t,\xi) + \int_0^t \partial_t^2\Phi_2(t,\tau,\xi)\hat{f}(\tau,\xi)\mathrm{d}\tau$$

from $(\partial_t\Phi_2)(t,t,\xi) = 1$. Hence if we denote $\mathrm{P} = \partial_t^2 + \frac{\mu}{1+t}\partial_t + |\xi|^2$ we have from $\mathrm{P}\Phi_2 = 0$
$$\mathrm{P} \int_0^t \Phi_2(t,\tau,\xi)\hat{f}(\tau,\xi)\mathrm{d}\tau$$
$$= \hat{f}(t,\xi) + \int_0^t \mathrm{P}\Phi_2(t,\tau,\xi)\hat{f}(\tau,\xi)\mathrm{d}\tau = \hat{f}(t,\xi).$$

Obviously function (4.4) satisfies zero initial conditions
$$\left[\int_0^t \Phi_2(t,\tau,\xi)\hat{f}(\tau,\xi)\mathrm{d}\tau\right]_{t=0} = 0,$$
$$\partial_t\left[\int_0^t \Phi_2(t,\tau,\xi)\hat{f}(\tau,\xi)\mathrm{d}\tau\right]_{t=0} =$$
$$\Phi_2(0,0,\xi)\hat{f}(t,\xi) + \left[\int_0^t \partial_t\Phi_2(t,\tau,\xi)\hat{f}(\tau,\xi)\mathrm{d}\tau\right]_{t=0} = 0.$$

Hence combining both representations by linearity we have

$$\hat{v}(t,\xi) = \sum_{i=1,2} \Phi_i(t,0,\xi)\hat{v}_i + \int_0^t \Phi_2(t,\tau,\xi)\hat{f}(\tau,\xi)d\tau. \qquad (4.5)$$

If we set

$$K_i(t,t_0,x) = (2\pi)^{-\frac{n}{2}} \mathcal{F}^{-1}_{\xi \to x}\big[\Phi_i(t,t_0,\xi)\big], \qquad i=1,2, \qquad (4.6)$$

Fourier transform in the sense of $\mathcal{S}'(\mathbb{R}^n)$, we have convoplution representation

$$v(t,x) = K_1(t,0,x) * v_1(x) + K_2(t,0,x) * v_2(x) + \int_0^t K_2(t,\tau,x) * f(\tau,x)d\tau. \qquad (4.7)$$

## 4.2 A-priori Estimates

We sketch only how to obtain a-priori estimates containing a nontrivial right hand side $f$. We restrict ourselves to the simpler case of estimate for the $L_q$-energy.

**Theorem 4.1.** *Let $u \in \mathcal{S}(\mathbb{R}^n)$. Let further $p \in (1,2]$ and $q$ such that $pq = p + q$. Then the estimate*

$$||(\partial_t, \nabla)K_2(t,\tau) * u||_q \lesssim \left(\frac{1+t}{1+\tau}\right)^{\max\{-\frac{n-1}{2}\left(\frac{1}{p}-\frac{1}{q}\right)-\frac{\mu}{2},\, -n\left(\frac{1}{p}-\frac{1}{q}\right)-1\}} ||u||_{L_{p,r}}$$

*holds for all $0 \leq \tau \leq t$ with $r = n\left(\frac{1}{p} - \frac{1}{q}\right)$.*

The key point is to include also the variable $\tau$ into the phase variables and in the definition of the zones. We use

$$Z_1 = \{(1+\tau)|\xi| \geq K\},$$
$$Z_2 = \{(1+t)|\xi| \geq K \geq (1+\tau)|\xi|\},$$
$$Z_3 = \{K \geq (1+t)|\xi|\},$$

and follow in each part the reasoning of Section 3.1.

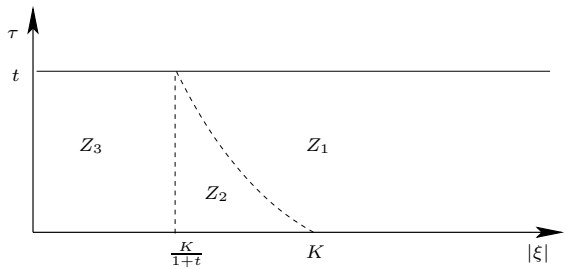

As consequence we obtain the following estimate of the $L_q$-energy of the solution to (4.1).

**Corollary 4.2.** *It holds*

$$\|(\partial_t, \nabla)v(t,\cdot)\|_q$$
$$\lesssim (1+t)^\delta \left[\|v_1\|_{L_{p,r+1}} + \|v_2\|_{L_{p,r}} + \int_0^t (1+\tau)^{-\delta}\|f(\tau,\cdot)\|_{L_{p,r}}\,d\tau\right] \quad (4.8)$$

*for $p \in (1,2]$, $q$ with $pq = p+q$, $r = n\left(\frac{1}{p} - \frac{1}{q}\right)$ and*

$$\delta = \max\left\{-\frac{n-1}{2}\left(\frac{1}{p}-\frac{1}{q}\right) - \frac{\mu}{2},\ -n\left(\frac{1}{p}-\frac{1}{q}\right) - 1\right\}.$$

*Remark* 4.1. The estimate is much weaker than the corresponding estimate for the free wave equation, which contains $(1+t-\tau)$ instead of $(1+t)/(1+\tau)$. The reason for this is the more complex structure of our model problem and the absence of the invariance under time translations.

# A  Appendix – Relation to Known Cases

## A.1  The case $\mu = 0$

The multiplier representation for the solution of (1.1) obtained in Section 2 can be simplified in the case $\mu = 0$. In this case we have $\rho = 1/2$ and hence the Bessel functions occuring in (2.21) reduce to trigonometric functions by

$$\mathcal{J}_{\frac{1}{2}}(\tau) = \sqrt{\frac{2}{\pi\tau}} \sin \tau \qquad (A.1)$$

$$\mathcal{J}_{-\frac{1}{2}}(\tau) = \sqrt{\frac{2}{\pi\tau}} \cos \tau. \qquad (A.2)$$

The structure of the multiplier simplifies by the addition theorem of trigonometric functions. This yields

$$\Phi_1(t,t_0,\xi) \qquad (A.3)$$
$$= \frac{\pi}{2}|\xi|\sqrt{(1+t)(1+t_0)}$$
$$\quad \left\{\mathcal{J}_{\frac{1}{2}}((1+t_0)|\xi|)\mathcal{J}_{\frac{1}{2}}((1+t)|\xi|) + \mathcal{J}_{-\frac{1}{2}}((1+t_0)|\xi|)\mathcal{J}_{-\frac{1}{2}}((1+t)|\xi|)\right\}$$
$$= \sin\big((1+t_0)|\xi|\big)\sin\big((1+t)|\xi|\big) + \cos\big((1+t_0)|\xi|\big)\cos\big((1+t)|\xi|\big)$$
$$= \cos\big((t-t_0)|\xi|\big). \qquad (A.4)$$

For the second multiplier one obtains in a similar way

$$\Phi_2(t, t_0, \xi) = \frac{\sin\left((t-t_0)|\xi|\right)}{|\xi|}. \tag{A.5}$$

These multipliers representing the solution of the wave equation are well-known.

The difference $t - t_0$ is natural in this represention because of the invariance of the differential operator under time translations. This can only be valid for $\mu = 0$.

For all even values of $\mu$ it is possible to represent the occuring Bessel functions by trigonometric functions, but only for $\mu = 2$ the representation simplifies significantly, for the reason cf. A.2.

In the general case an addition theorem for Bessel functions like that for trigonometric ones is not available.

## A.2 Transformation to a Klein-Gordon Equation with Variable Mass

It is possible to transform the equation

$$\Box v + \frac{\mu}{1+t} v_t = 0 \tag{A.6}$$

to a Klein-Gordon type equation. We follow [Mat76] and [Rei00]. We have

$$\exp\left\{-\frac{1}{2}\int_0^t \frac{\mu}{1+\tau} d\tau\right\} = (1+t)^{-\frac{\mu}{2}}. \tag{A.7}$$

If we set

$$v = (1+t)^{-\frac{\mu}{2}} w \tag{A.8}$$

we obtain

$$\begin{aligned}
0 &= \left((1+t)^{-\frac{\mu}{2}} w\right)_{tt} - (1+t)^{-\frac{\mu}{2}} \Delta w + \tfrac{\mu}{1+t}\left((1+t)^{-\frac{\mu}{2}} w\right)_t \\
&= \left(-\tfrac{\mu}{2}(1+t)^{-\frac{\mu}{2}-1} w + (1+t)^{-\frac{\mu}{2}} w_t\right)_t - (1+t)^{-\frac{\mu}{2}} \Delta w \\
&\quad + \tfrac{\mu}{1+t}\left(-\tfrac{\mu}{2}(1+t)^{-\frac{\mu}{2}-1} w + (1+t)^{-\frac{\mu}{2}} w_t\right) \\
&= -\tfrac{\mu}{2}\left(-\tfrac{\mu}{2}-1\right)(1+t)^{-\frac{\mu}{2}-2} w - \mu(1+t)^{-\frac{\mu}{2}-1} w_t + (1+t)^{-\frac{\mu}{2}} w_{tt} \\
&\quad - (1+t)^{-\frac{\mu}{2}} \Delta w - \tfrac{\mu^2}{2}(1+t)^{-\frac{\mu}{2}-2} w + \mu(1+t)^{-\frac{\mu}{2}-1} w_t \\
&= (1+t)^{-\frac{\mu}{2}} \left\{\Box w + \tfrac{\mu(2-\mu)}{4}(1+t)^{-2} w\right\}.
\end{aligned}$$

After multiplication with $(1+t)^{\frac{\mu}{2}}$ we obtain for $w$ the Klein-Gordon equation

$$\Box w + \tfrac{\mu(2-\mu)}{4}(1+t)^{-2}w = 0 \tag{A.9}$$

with mass tending to zero as $t$ tends to infinity. The sign of the mass term depends on the size of $\mu$. For $\mu \in (0,2)$ the mass is positive and for $\mu > 2$ negative. For the special case $\mu = 2$ we obtain a wave equation.

*Remark* A.1. This gives an interpretation for the change in the behaviour of the energy estimate for $\mu = 2$. In the forthcoming papers we will develop this in more detail for more general classes of weak dissipation.

Hence if $\mu = 2$ known a-priori estimates for the wave equation (cf. [Str70] or [vW71]) lead immediately to

$$||v(t,\cdot)||_\infty \leq C_\epsilon (1+t)^{-\frac{n-1}{2}-1}\big(||v_1||_{W_1^{n+\epsilon}} + ||v_2||_{W_1^{n-1+\epsilon}}\big), \tag{A.10a}$$

$$||v(t,\cdot)||_2 \leq C\big(||v_1||_2 + ||v_2||_{H^{-1}}\big), \tag{A.10b}$$

and for the derivatives

$$||(v_t(t,\cdot),\nabla v(t,\cdot))||_\infty \leq C_\epsilon (1+t)^{-\frac{n-1}{2}-1}\big(||v_1||_{W_1^{n+1+\epsilon}} + ||v_2||_{W_1^{n+\epsilon}}\big), \tag{A.11a}$$

$$||(v_t(t,\cdot),\nabla v(t,\cdot))||_2 \leq C(1+t)^{-1}\big(||v_1||_{H^1} + ||v_2||_2\big). \tag{A.11b}$$

These estimates are related to the estimates from Section 3.2 and 3.3. The occurence of the parameter $\epsilon$ in these $L_1$–$L_\infty$ estimates is necessary. The space $W_1^n(\mathbb{R}^n)$ itself is not embedded in $L_\infty(\mathbb{R}^n)$.

## Acknowledgements

The author has to thank the Institute of Mathematics of the Tsukuba University, especially Prof. K. Kajitani, for the invitation and the hospitality during his stay in February 2002 as well as Prof. M. Reissig for many discussions and helpful remarks during these calculations and the preparation of this paper.